\documentclass{amsart}

\usepackage{graphicx}
\usepackage[centertags]{amsmath}
\usepackage{amsfonts}
\usepackage{amssymb}
\usepackage{amsthm}
\usepackage{newlfont}
\usepackage{hyperref}

\newtheorem{thm}{Theorem}
\newtheorem{cor}{Corollary}

\newtheorem{prop}{Proposition}

\theoremstyle{definition}

\theoremstyle{remark}

\newcommand{\DM}{\mathrm{F}}
\newcommand{\Alt}{\mathrm{Alt}}
\newcommand{\N}{\mathbb{N}}
\newcommand{\Z}{\mathbb{Z}}

\begin{document}

\title{Arbitrarily large residual finiteness growth}
\author{Khalid Bou-Rabee}
\address{Khalid Bou-Rabee, Department of Mathematics, The University of Michigan, 
2074 East Hall, Ann Arbor, MI 48109-1043}
\email{khalidb@umich.edu}

\author{Brandon Seward}
\address{Brandon Seward, Department of Mathematics, The University of Michigan, 
2074 East Hall, Ann Arbor, MI 48109-1043}
\email{bseward@umich.edu}

\keywords{residual finiteness growth, residually finite, B. H. Neumann groups.}
\subjclass{20E26}

\begin{abstract}
The residual finiteness growth of a group quantifies how well approximated the group is by its finite quotients. In this paper, we construct groups with arbitrarily large residual finiteness growth.
We also demonstrate a new relationship between residual finiteness growth and some decision problems in groups, which we apply to our new groups.
\end{abstract}
\maketitle

Let $\Gamma$ be a residually finite group generated by a finite set $X$.
An element $\gamma \in \Gamma$ is \emph{detected} by a finite quotient $Q$ of $\Gamma$ if there exists a homomorphism $\phi : \Gamma \to Q$ such that $\gamma \notin \ker\phi$.
Set $\DM_{\Gamma, X}(n)$ to be the minimal natural number $N$ such that any element of  word length $\leq n$ with respect to $X$ can be detected in a quotient $Q$ of cardinality $\leq N$.
The study of $\DM_{\Gamma, X}$, the \emph{residual finiteness growth}, is an area of active research (see, for instance, \cite{B10}, \cite{B11},  \cite{BM1},  \cite{BK10}, and \cite{KM11}).
Our main result is the following

\begin{thm}
\label{MainTheorem}
For any function $f : \N \rightarrow \N$ there exists a residually finite group $B_f$ and a two element generating set $X$ for $B_f$ such that $\DM_{B_f, X}(n) \geq f(n)$ for all $n \geq 8$.
\end{thm}
Our proof is constructive, elementary, and demonstrates the versatility of B. H. Neumann's construction of a continuous family of non-isomorphic two generated groups \cite{N37}.
Our work is of the same flavor as \cite{LPS96} and \cite{P04}, where Neumann's construction is modified to find groups with interesting subgroup growth properties.

Our remaining results exhibit how residual finiteness growth may be applied to illicit decision problem data on groups. 
The first of these shows that a variant of Theorem \ref{MainTheorem} does not exist for groups that are finitely presented.
In particular, the main result on residual finiteness growth stated in \cite{KMS13} is the strongest of its kind for finitely presented groups.

\begin{prop} \label{DecTheorem1}
If $\Gamma$ is a residually finite finitely presented group, then the residual finiteness growth function for $\Gamma$ is computable.
\end{prop}

O.~Kharlampovich showed that there exists a finitely presented solvable group with undecidable word problem \cite{K81} of which many variations have been found (see \cite{B83}).
In \cite[IV.4.3 Theorem 19]{LS77} it is shown that $F_2 \times F_2$ does not have decidable membership problem with respect to finitely generated subgroups. 
We address a weaker membership problem for the groups $B_f$.
For a finitely generated group $G$, we say the \emph{finite-index normal subgroup membership problem in $G$ is solvable} if there exists an algorithm which for every element $g \in G$ and every finite-index normal subgroup $H \lhd G$ will determine in in a finite amount of time whether $g \in H$.

\begin{thm} \label{DecTheorem2}
If the finite-index normal subgroup membership problem in $\Gamma$ is solvable and the word problem for $\Gamma$ is solvable, then the residual finiteness growth function for $\Gamma$ is computable.
\end{thm}

We remark that recursively presented finitely generated residually finite groups may have undecidable word problem (see \cite{M74} or \cite{D74}). In \cite{BM09} Baumslag and Miller studied finite presentability and the solvability of the word problem in Neumann's groups. Our groups are a variant of Neumann's groups and thus not formally covered by the findings in \cite{BM09}. However, by applying Theorem \ref{DecTheorem2} we are able to use the residual finiteness growth function to obtain information on finite presentability and decision problems for the groups we construct. We believe this is the first time residual finiteness growth has been used in this manner.

\begin{cor}
If $f : \N \rightarrow \N$ is a function which grows faster than every computable function, then the group $B_f$:
\begin{enumerate}
\item is not finitely presented 
\item has unsolvable finite-index normal subgroup membership problem or unsolvable word problem.
\end{enumerate}
\end{cor}

\section{Acknowledgements}
The first author is grateful to Benson Farb and Ben McReynolds for their endless support and great ideas.
The first author was partially supported by NSF RTG grant DMS-0602191.
The second author would like to thank Jay Williams for an insightful discussion of Neumann's construction. The second author was supported by a NSF Graduate Student Research Fellowship under Grant No. DGE 0718128.

\section{Proof of Theorem \ref{MainTheorem}}

\begin{proof}
Loosely speaking, the key idea of this proof is to construct a finitely generated infinite group containing infinitely many alternating groups that do not wander far from the identity element. Hence, the structure of the finite alternating groups $\Alt(n)$ for $n$ odd play a significant role. Recall that for odd $n \geq 5$ the group $\Alt(n)$ is simple and is generated by the following two elements \cite[III.B.26]{H00}:
$$\alpha_n = (1, 2, 3, \ldots, n)$$
$$\beta = (1, 2, 3).$$
Observe that
$$\alpha_n^m \beta \alpha_n^{-m} = (m+1, m+2, m+3),$$
where the coordinates are understood to be modulo $n$. Therefore $\alpha_n^m \beta \alpha_n^{-m}$ and $\beta$ commute precisely when $m \not\equiv \pm 1, \pm 2 \mod n$. So
$$[\alpha_n^m \beta \alpha_n^{-m}, \beta] = 1 \Longleftrightarrow m \not\equiv \pm 1, \pm 2 \mod n.$$
By playing off of congruence properties of integers, we will build a group in which one can arrive at arbitrarily large alternating groups with words of relatively short length. In the next paragraph we construct the integral sequences which will serve this purpose.

Without loss of generality, we may suppose that $f : \N \rightarrow \N$ is increasing and $f(n) \geq n$ for all $n \in \N$. We begin by constructing three functions $p, q, d: \N \rightarrow \N$ satisfying the following conditions for all $k \in \N$:
\begin{enumerate}
\item[\rm (i)] $p: \N \rightarrow \N$ is strictly increasing, takes only odd values, and $p(1) = 1$;
\item[\rm (ii)] $q(k) > 2$ is odd;
\item[\rm (iii)] $d(k) = p(k) q(k) + 2$;
\item[\rm (iv)] $d(k)$ is an odd integer greater than $2 \cdot f(4 p(k+1) + 4)$;
\item[\rm (v)] for all $i \leq k$, $p(k + 1) q(i) \not\equiv \pm 1, \pm 2 \mod d(i)$.
\end{enumerate}
To begin we set $p(1) = 1$, $p(2) = 3$, $d(1) = 2 \cdot f(16) + 1$, and $q(1) = 2 \cdot f(16) - 1$. Note that $q(1) \equiv -2 \mod d(1)$ and hence
$$p(2) q(1) \equiv 3 q(1) \equiv -6 \not\equiv \pm 1, \pm 2 \mod d(1)$$
since $d(1) \geq 33$. Now suppose that $q(k-1)$, $d(k-1)$, and $p(k)$ have been defined. Let $\ell$ be the least common multiple of $d(1), d(2), \ldots, d(k-1)$. Set
$$q(k) = 2 \cdot f(4 p(k) + 40 \ell + 4) + 1$$
and set $d(k) = p(k) q(k) + 2$. Then $d(k)$ is odd since both $p(k)$ and $q(k)$ are odd. Also, $q(k)$ and $d(k)$ must be relatively prime since they are both odd and any number dividing both $d(k)$ and $q(k)$ must divide $2$. So $q(k)$ is not a zero divisor in $\Z / d(k) \cdot \Z$ and hence there are precisely four values of $m$ modulo $d(k)$ satisfying the equation
$$m \cdot q(k) \equiv \pm 1, \pm 2 \mod d(k).$$
Therefore we can pick $p(k+1)$ out of the five odd numbers $\{p(k) + 2 \ell, p(k) + 4 \ell, \ldots, p(k) + 10 \ell\}$ so that $p(k+1) q(k) \not\equiv \pm 1, \pm 2 \mod d(k)$ (we use here the fact that these five numbers are distinct modulo $d(k)$ as $p(k) + 10 \ell < d(k)$). With these definitions clauses (i), (ii), (iii), and (iv) continue to be satisfied. Clause (v) is immediate from our definition when $i = k$. For $i < k$ we have $p(k+1) \equiv p(k) \mod d(i)$ and hence
$$p(k+1) q(i) \equiv p(k) q(i) \not\equiv \pm 1, \pm 2 \mod d(i).$$
This completes the inductive definitions of $p$, $q$, and $d$.

Now we construct the group $B_f$. Let
$$G = \prod_{k \in \N} \Alt(d(k)).$$
We write $G_k$ for the canonical copy of $\Alt(d(k))$ in $G$, and we let $\pi_k : G \rightarrow \Alt(d(k))$ be the canonical homomorphism. Set
$$s = (s_1, s_2, s_3, \ldots) \in G$$
where $s_k = \alpha_{d(k)}^{q(k)}$ and also set
$$t = (\beta, \beta, \beta, \ldots) \in G.$$
We set $B_f = \langle s, t \rangle$ and $X = \{s, t\}$. Clearly $B_f$ is $2$-generated as required. By restricting the homomorphisms $\pi_k$, $k \in \N$, we immediately see that $B_f$ is residually finite.

As stated previously, the main idea is that one can enter into large finite simple subgroups by using relatively short words. The precise claim is the following. For every $k \in \N$, $G_k$ is a subgroup of $B_f$ and contains an element of $X$-word-length $4 p(k) + 4$. Fix $k \in \N$ and consider the word
$$w(u, v) = [u^{p(k)} v u^{-p(k)}, v].$$
Clearly $w(s, t)$ has $X$-word-length $4 p(k) + 4$. For $i < k$ we have $p(k) q(i) \not\equiv \pm 1, \pm 2 \mod d(i)$ by clause (v) and hence
$$\forall 1 \leq i < k \quad \pi_i(w(s, t)) = w(s_i, \beta) = 1.$$
For $i > k$ we have
$$2 < q(i) < p(k) q(i) < p(i) q(i) = d(i) - 2$$
by clauses (i), (ii), and (iii). Hence $p(k) q(i) \not\equiv \pm 1, \pm 2 \mod d(i)$ and thus
$$\forall i > k \quad \pi_i(w(s, t)) = w(s_i, \beta) = 1.$$
However, by clause (iii) we have $p(k) q(k) \equiv -2 \mod d(k)$ and thus
$$\pi_k(w(s, t)) = w(s_k , \beta) \neq 1.$$
Therefore $w(s, t) \in G_k$. This proves half of our claim. So now we have $B_f \cap G_k \neq \{1\}$. As $G_k$ is simple ($d(k) \geq 5$ by clauses (ii) and (iii)), it suffices to show that $B_f \cap G_k$ is a normal subgroup of $G_k$. Fix $\gamma \in B_f \cap G_k$ and $g \in G_k$. We will show that $g \gamma g^{-1} \in B_f \cap G_k$. As stated at the beginning of this proof, $\Alt(d(k))$ is generated by $\alpha_{d(k)}$ and $\beta$. As observed earlier, clauses (ii), (iii), and (iv) imply that $q(k)$ and $d(k)$ are relatively prime. Thus $\langle s_k \rangle = \langle \alpha_{d(k)} \rangle$ and hence $\Alt(d(k)) = \langle s_k, \beta \rangle$. So $\pi_k(g) = w'(s_k, \beta) \in \langle s_k, \beta \rangle$ for some word $w'(u, v)$. Setting $\lambda = w'(s, t) \in B_f$ we have
$$\pi_k(\lambda) = w'(s_k, \beta) = \pi_k(g).$$
As $\gamma, g, g \gamma g^{-1} \in G_k$ we have $\pi_i(\gamma) = 1$ and $\pi_i(g \gamma g^{-1}) = 1$ for $i \neq k$. So we have
$$\pi_k(\lambda \gamma \lambda^{-1}) = \pi_k(\lambda) \pi_k(\gamma) \pi_k(\lambda)^{-1} = \pi_k(g) \pi_k(\gamma) \pi_k(g)^{-1} = \pi_k(g \gamma g^{-1})$$
and for $i \neq k$
$$\pi_i(\lambda \gamma \lambda^{-1}) = \pi_i(\lambda) \pi_i(\gamma) \pi_i(\lambda)^{-1} = \pi_i(\lambda) \cdot 1 \cdot \pi_i(\lambda)^{-1} = 1 = \pi_i(g \lambda g^{-1}).$$
Therefore $g \gamma g^{-1} = \lambda \gamma \lambda^{-1} \in B_f$. Thus $g \gamma g^{-1} \in B_f \cap G_k$ and by the simplicity of $G_k$ we conclude that $B_f \cap G_k = G_k$. This proves the claim.

Now we finish the proof by showing that $\DM_{B_f, X}(n) \geq f(n)$ for all $n \geq 8$. So fix $n \geq 8$. Since $4 p(1) + 4 = 8$ and $p: \N \rightarrow \N$ is strictly increasing, we can find $k \in \N$ with $4 p(k) + 4 \leq n < 4 p(k + 1) + 4$. By the previous paragraph, there is $\gamma \in B_f$ of $X$-word-length $4 p(k) + 4 \leq n$ with $\gamma \in G_k \leq B_f$. Since $G_k$ is simple, any normal subgroup of $B_f$ not containing $\gamma$ must meet $G_k$ trivially. This implies that distinct members of $G_k$ would lie in distinct cosets of such a normal subgroup. Therefore by clause (iv)
$$\DM_{B_f, X}(n) \geq |G_k| = \frac{1}{2} d(k)! > \frac{1}{2} d(k) \geq f(4 p(k+1) + 4) > f(n).$$
\end{proof}

\section{Proof of Proposition \ref{DecTheorem1}}

\begin{proof}
Fix a finite generating set $X$ for $\Gamma$. We give an algorithm for computing $\DM_{\Gamma, X}$.
Let $n \in \N$ be given.
For each such $n$, we enumerate all the nontrivial elements $\{ \gamma_k \}$ in $B_{\Gamma,X}(n)$ (we can do this because a finitely presented residually finite group has solvable word problem).
So long as there exists an algorithm listing all maps $\Gamma \to Q$, $Q$ finite, we may enumerate all such maps $\phi_1, \phi_2, \ldots$ such that their images are nondecreasing (done in the next paragraph).
For each $\phi_i$ we can decide whether $\phi_i(\gamma_k)$ is trivial by following the word determining $\gamma_k$ in the Cayley graph of $Q$ (the Cayley graph of any finite group in terms of a fixed generating set is given by an algorithm that terminates in finite time).
Thus, we can determine the least $m$ so that for every non-trivial $\gamma_k \in B_{\Gamma, X}(n)$ there is $i \leq m$ with $\phi_i(\gamma_k) \neq 1$. So $\DM_{\Gamma,X}(n)$ is computable as it is the order of the image of $\phi_m$.

To finish, we need only show that there is an algorithm listing all maps from $\Gamma$ onto finite groups. Let $\Gamma = \langle X | R \rangle$ be a finite presentation of $\Gamma$. Fix an ordering for $X$ and set $k = |X|$.
Enumerate all finite groups generated by $k$ (or fewer) elements by $Q_i$ and enumerate all ordered generating sets of $Q_i$ of cardinality $k$ by $Y_{ij}$.
For each $Y_{ij}$ the unique order-preserving bijection $\phi: X \rightarrow Y_{ij}$ extends to a homomorphism $\phi: \Gamma \rightarrow Q_i$ if and only if $\phi(r) = 1$ for every $r \in R$.
Since both $R$ and $Y_{ij}$ are finite this property can be checked by an algorithm in finite time.
This gives us an enumeration of all the finite quotients of $\Gamma$ in increasing order.
\end{proof}

\section{Proof of Theorem \ref{DecTheorem2}}

\begin{proof}
Since we are assuming $\Gamma$ has solvable word problem, we may, as in the proof of Theorem \ref{DecTheorem1}, only show that there exists an algorithm that enumerates all maps onto finite quotients of $\Gamma$.
Again, let $X$ be an ordered finite generating set for $\Gamma$ and set $k = |X|$.
As in the previous proof, enumerate all finite groups $Q_i$ (in increasing order) which are generated by $k$ (or fewer) elements. For each such group $Q_i$, let $Y_{ij}$ enumerate all ordered generating sets of cardinality $k$.
Again, for each $Y_{ij}$ we consider the unique order-preserving bijection $\phi_{ij}: X \rightarrow Y$ and seek to determine if $\phi_{ij}$ extends to a homomorphism $\phi_{ij}: \Gamma \rightarrow Q_i$.
Let $F_X$ be the free group over the alphabet $X$. Clearly $\phi_{ij}$ naturally extends to a homomorphism $\phi_{ij}: F_X \rightarrow Q_i$. Since $Q_i$ is finite, one can algorithmically find a finite set $R_{ij} \subseteq F_X$ with $\langle R_{ij} \rangle = \ker \phi_{ij}$ (equivalently, $Q_i$ has presentation $Q_i = \langle Y_{ij} | R_{ij} \rangle$). Let $\pi : F_X \rightarrow \Gamma$ be the obvious homomorphism and let $K = \ker \pi$. Clearly $\phi_{ij}$ extends to a homomorphism $\phi_{ij} : \Gamma \rightarrow Q_i$ if and only if $K \leq \langle R_{ij} \rangle$. So it suffices to determine if $K \leq \langle R_{ij} \rangle$. We have
$$K \leq \langle R_{ij} \rangle \Longleftrightarrow |F_X : K \cdot \langle R_{ij} \rangle| = |F_X : \langle R_{ij} \rangle| \Longleftrightarrow |\Gamma : \langle \pi(R_{ij}) \rangle| = |F_X : \langle R_{ij} \rangle|.$$
By using preimages of elements of $Q_i$, we can compute a set of coset representatives $T_{ij} \subseteq F_X$ for the subgroup $\langle R_{ij} \rangle$. Now we can determine the index of $\langle \pi(R_{ij}) \rangle$ in $\Gamma$ by counting the number of $t \in T_{ij}$ with $\pi(t) \in \langle \pi(R_{ij}) \rangle$. Since $\Gamma$ has solvable finite-index normal subgroup membership problem and $T_{ij}$ is finite, there is an algorithm which can do this in finite time. Thus we can enumerate the finite quotients of $\Gamma$. The remainder of the argument is identical to the first paragraph of the proof of Proposition \ref{DecTheorem1}.

\end{proof}

\end{document}